\documentclass[10pt, twoside, a4paper]{article}
\usepackage[utf8]{inputenc}
\usepackage[T1]{fontenc}
\usepackage{amsmath, amssymb, amsthm} 
\usepackage{amstext, amsfonts, a4}
\usepackage[english]{babel}
\usepackage{stmaryrd}
\usepackage{bbm}
\usepackage{hyperref}
\usepackage[margin=3cm]{geometry}

\theoremstyle{plain}
	\newtheorem{theorem}{Theorem}[section]
		
	\newtheorem{lemma}{Lemma}[section]
	
\theoremstyle{definition}
	\newtheorem{definition}{Definition}[section]

\theoremstyle{remark}

\DeclareMathOperator{\supp}{supp}

       \def\HH{{\mathbb H}}          \def\RR{{\mathbb R}}        \def\Z\int{{\mathbb Z}} 

    \def\cB{{\mathcal B}}  \def\cN{{\mathcal N}}  \def\cC{{\mathcal C}}  \def\cO{{\mathcal O}} \def\cU{{\mathcal U}}         \def\cF{{\mathcal F}} \def\cL{{\mathcal L}}    

\begin{document}

\title{Smooth traveling-wave solutions to the inviscid
surface quasi-geostrophic equations}
\author{
\renewcommand{\thefootnote}{\arabic{footnote}}
Ludovic Godard-Cadillac\footnotemark[1]}
\footnotetext[1]{Sorbonne Universit\'e, Laboratoire Jacques-Louis Lions, 4, Place Jussieu, 75005 Paris, France.\\ E-mail: {\tt ludovic.godard-cadillac@sorbonne-universite.fr}}
\date{\today}
\maketitle

\begin{abstract}
In a recent article by Gravejat and Smets~\cite{Gravejat_Smets_2019}, it is built smooth solutions to the inviscid
surface quasi-geostrophic equation that have the form of a traveling wave. In this article we
work back on their construction to provide solution to a more general class of quasi-geostrophic
equation where the half-laplacian is replaced by any fractional laplacian.
\end{abstract}

\tableofcontents

\section{Presentiation of the problem}
\subsection{The quasi-geostrophic equations}
We consider the general transport equation for the vorticity of an incompressible fluid in dimension
2.
\begin{equation}\label{eq:SQG transport}
\frac{\partial\theta}{\partial t} + v \cdot \nabla\theta = 0,
\end{equation}
where $\theta : \RR^2\times\RR_+ \to \RR$ is the called the \emph{active scalar} and $v : \RR^2\times\RR_+ \to \RR$
is the \emph{velocity} of the fluid This equation tells that the active scalar is transported by the induced velocity. Since this
velocity $v$ is divergence free (incompressibility condition), it is convenient to relate $v$ and $\theta$ through
a stream function $\psi : \RR^2 \times \RR_+ \to \RR.$ The generalized inviscid surface quasi-geostrophic equation
corresponds to a stream function that verifies
\begin{equation}\label{eq:SQG stream} v = \nabla^\perp\psi\qquad \text{and}\qquad(-\Delta)^s\psi = \theta, \end{equation}
with $s \in]0, 1[$ and $\perp$ denotes the rotation in the plane of angle $\frac{\pi}{2}$. The equations~\eqref{eq:SQG transport}-\eqref{eq:SQG stream} are the generalized inviscid surface quasi-geostrophic
equations. In the particular case $s \to 1$ we obtain the well-known 2D Euler equation written in
term of vorticity and stream function. Another important case is $s =\frac{1}{2}$ which correspond to
the work made in~\cite{Gravejat_Smets_2019} that we generalize here. This case is the standard surface quasi-geostrofic
equation that first appeared as a limit model in the context of geophysical flows~\cite{Pedlowsky_1987}\cite{Vallis_2006}. 
These equations are used to model a fluid in a rotation frame with stratified density and velocity and submitted to Brunt-Väisälä thermal oscillations. This models leads to~\eqref{eq:SQG transport}-\eqref{eq:SQG stream} using the Cafferelli-Silverstre theory on fractional laplace operator~\cite{Caffarelli_Silvestre_2007}.
The case of the exponent $s=\frac{1}{2}$ corresponds to the case of a Brunt-Väisälä frequency $N$ that does not depend on the height. Other exponents for the fractional Laplace operator corresponds to different profiles for $N$. These equations has been
intensely investigated since the work of Constantin, Majda and Tabak~\cite{Constantin_Majda_Tabak_1994} on the case $s=\frac{1}{2}$ where they pointed
out the mathematical links that arises between (SQG-$\frac{1}{2}$) and the Euler
equation in dimension 3. Besides stationary solution, given by a radially symmetric rearrangement on the
active scalar, the only two known examples of global smooth solutions where built by Castro,
Córdoba and Gómez-Serrano~\cite{Castro_Cordoba_Gomez_2016} on the one hand and by Gravejat and Smets~\cite{Gravejat_Smets_2019} on the other hand
with two different techniques. The article of Castro, Córdoba and Gómez-Serrano also provides a
wide bibliography related on SQG and its Cauchy problem.
In this work we generalize the result and the construction provided by~\cite{Gravejat_Smets_2019} to the more general
equations~\eqref{eq:SQG transport}-\eqref{eq:SQG stream}  with a fixed $s \in]0, 1[$. The idea consists in looking for solutions that have the
form of traveling waves with a positive speed $c$ in direction $z$. In short, solutions of
the form
\begin{equation}
\theta (r, z, t) = \Theta (r, z - ct),\qquad v (r, z, t) = V (r, z - ct),\qquad
\psi (r, z, t) = \Psi (r, z - ct).
\end{equation}
We inject this form of solution in~\eqref{eq:SQG transport}-\eqref{eq:SQG stream} and we get
\begin{equation}\begin{split}
0&=\frac{\partial}{\partial t} (\Theta (r, z - ct)) + V(r, z - ct).\nabla\Theta (r, z - ct)\\
&= -ce_z.\nabla\Theta (r, z - ct) + V(r, z - ct).\nabla\Theta (r, z - ct)\\
&= -ce_z.\nabla\Theta (r, z - ct) + \nabla^\perp\Psi (r, z - ct).\nabla\Theta (r, z - ct),
\end{split}
\end{equation}
where $(e_r,e_z)$ denotes the canonical basis of $\RR^2$.
This leads to the orthogonality condition
\begin{equation}\label{eq:orthobonality condition}
\Big(\nabla\psi - c e_r\Big)^\perp
\cdot\nabla\Theta = 0,
\end{equation}
with the remark that $e_r^\perp = e_z$. In other words, the two vectors $\nabla\Theta$ and $\nabla\Psi-c e_r$ 
must be collinear. Following an idea from Arnold~\cite{Arnold_1978}, Condition~\eqref{eq:orthobonality condition} is immediately verified if $\Theta$ has the form
\begin{equation}
\Theta (r, z) = f\Big(\Psi (r, z) - cr - k\Big)
\end{equation}
because in this case
\begin{equation}
\nabla\Theta(r, z) = f'\big(\Psi (r, z) - cr - k\big)\cdot\big(\nabla\Psi (r, z) - c e_r\big)
\end{equation}
which gives~\eqref{eq:orthobonality condition}. We now make consider the ansatz of a symmetry relatively to the $z$-axis that takes the
form
\begin{equation}\label{eq:symmetry hypothesis} \Psi (-r, z) = -\Psi (r, z).
\end{equation}
This implies that $\Theta(-r,z)=-\Theta(r,z)$ and if we denote $V=(V_1,V_2)$ the two components of the velocity profile, then $V_1(-r,z)=-V_1(r,z)$ and $V_2(-r,z)=V_2(r,z)$.
More precisely, we impose the following ansatz
\begin{equation}\label{eq:ansatz}
\Theta (r, z) = \left\{\begin{split}&f (\Psi (r, z) - cr - k)\qquad\text{if } r \geq 0,\\
-&f (-\Psi (r, z) + cr - k)\qquad\text{otherwise} 
\end{split}\right.
\end{equation}
where $f : \RR \to \RR$ is a smooth function supported in $\RR_+$ (to avoid a singularity at $x = 0$) with the
condition $k > 0$. Using the stream equations~\eqref{eq:SQG stream} we obtain
\begin{equation}\label{eq:frac lap semi-lin}
(-\Delta)^s\Psi\;(r,z) = \left\{\begin{split}&f (\Psi (r, z) - cr - k)\qquad\text{if } r \geq 0,\\
-&f (-\Psi (r, z) + cr - k)\qquad\text{otherwise} 
\end{split}\right.
\end{equation}

\subsection{Variational formulation}
The studied equation is variational and its solutions are the critical points of
\begin{equation}
E(\Psi):= \frac{1}{2}\int_{\RR^2}\Psi (-\Delta)^s \Psi -
\int_{\HH} F (\Psi - cr - k) + \int_{\HH^c} F (\Psi + cr - k),
\end{equation}
where $\HH:=\{x=(r, z) \in\RR^2 : r \geq 0\}$ and $F(\xi):= \int_0^\xi
f (\xi') d\xi'$.
We are going to build a critical point
of E using the technique of the Nehari manifold (defined later). For that purpose, since the
choice of $f$ is free, we are imposing on this function the following properties
\begin{align}
&\bullet\quad f \in\cC^\infty (\RR, \RR),\qquad f_{|\RR_-} = 0 \quad\text{and}\quad f_{|\RR^\ast_+}> 0, \label{eq:hypothesis 1}\\
&\bullet\quad \exists\; \nu \in\bigg]1,\;\frac{1+s}{1-s}\bigg[,\; \forall \; \xi \geq 0,\quad f(\xi) \leq C \xi^\nu, \label{eq:hypothesis 2}\\
&\bullet\quad \exists\; \mu \in ]1, \nu[ ,\; \forall \; \xi \geq 0,\quad \mu\,f (\xi) \leq \xi f'(\xi). \label{eq:hypothesis 3}
\end{align}
This last hypothesis on the variations of $f$ is equivalent to the hypothesis that the function
\begin{equation}
\xi \longmapsto \frac{f(\xi)}{\xi^\mu}.
\end{equation}
is non-decreasing on $\RR_+$.
In particular and since $\mu > 1$,
\begin{equation}
\forall \; \xi_0 \geq 0,\quad \xi \in \RR_+ \longmapsto\frac{f(\xi)}{\xi + \xi_0}
\end{equation}
is increasing and diverging at infinity.
Examples of functions that satisfies these three hypothesis~\eqref{eq:hypothesis 1}\eqref{eq:hypothesis 2}\eqref{eq:hypothesis 3} are the functions
\begin{equation}
\xi \longmapsto \xi^\nu\,e^{-\frac{1}{\xi}}\,\mathbbm{1}_{\RR_+}(\xi),
\end{equation}
with $\nu \in [\mu, \nu]$. Given the hypothesis \eqref{eq:hypothesis 1} and \eqref{eq:hypothesis 2}, the functional $E$ is well-defined on the Hilbert
space
\begin{equation}
X^s:= L^\frac{2}{1-s} \cap \dot{H}^s(\RR^2)
\end{equation}
with the scalar product induced by $\dot{H}^s$ given by
\begin{equation}
\left<\Phi, \Psi\right>_{X^s} := p.v. \int_{\RR^2}\int_{\RR^2}\frac{\big(
\Phi (x) - \Phi (y)
\big)\big(\Psi (x) - \Psi (y)
\big)}{
|x - y|^{2(1+s)}}\,
dx\,dy
\end{equation}
where $p.v.$ refers to the principal value of the singularity of the kernel $(x,y)\mapsto1/|x - y|^{2(1+s)}.$ For further work, we make use of the notations $x=(r_x,z_x)$ and $y=(r_y,z_y)$ to distinguish the coordinates of $x$ and the coordinates of $y$.
We recall here that the Gagliardo half-norms defining the spaces $\dot{W}^{s,p}$ are given in
general by
\begin{equation}
|\Phi|_{W^{s,p}}^p:= p.v. \int_{\RR^d}\int_{\RR^d}\frac{\big|\Phi(x)-\Phi(y)\big|^p}{|x - y|^{d+sp}}\,dr \,dz. \label{eq:gagliardo half-norms}
\end{equation}
For the rest of the work we refer $E$ as being the ``\emph{energy}'' of the problem although this energy does
not correspond to a physical energy. We remark that it is invariant by the action of the group of symmetry generated by~\eqref{eq:symmetry hypothesis}.
 We denote by $X^s_{sym}$ the subspace of $X^s$ made with the functions
that are left invariant by the action of this symmetry group.
\begin{equation}
X^s_{sym} :=\{\Psi \in X^s\;:\;\forall \; (r, z) \in \RR^2, \Psi (-r, z) = -\Psi (r, z)\}.
\end{equation}
It follows from the Palais principle of symmetric criticality~\cite{Palais_1979} that any critical point of $E$ on $X^s$ actual belongs to $X^s_{sym}$. We can therefore restrict our investigations to the subspace $X^s_{sym}$, inside which the energy can be rewritten
\begin{equation}
E (\Psi) = \frac{1}{2}\|\Psi\|_{X^s}^2-2 V (\Psi) 
\end{equation}
with
\begin{equation}
V(\Psi):=\int_{\HH}F(\Psi - cr - k).
\end{equation}

\subsection{Nehari Manifold and presentation of the main result}

The Nehari manifold associated to the energy E is defined by
\begin{equation}
\cN =\{\Psi \in X^s_{sym}\setminus\{0\}\;:\;E'(\Psi) (\Psi) = 0\},
\end{equation}
so that $\Psi \in \cN$ implies
\begin{equation}
\int_{\RR^2}\Psi (-\Delta)^s \Psi - 2\int_{\HH} f (\Psi - cr - k) \Psi = 0. \label{eq:energy on Nehari}
\end{equation}
It is proven after that the Nehari manifold $\cN$ is a sub-manifold of $X^s_{sym}$ non empty, of regularity
$\cC^1$ without boundary. The main result of this article is the following theorem.
\begin{theorem}\label{thrm:}
Let $c$ and $k$ positive. Let $f : \RR \to\RR$ verifying \eqref{eq:hypothesis 1}, \eqref{eq:hypothesis 2} and \eqref{eq:hypothesis 3}. 

Then the energy
E admits a minimizer $\Psi \neq 0$ on $\cN$ . As a consequence there exist a non-trivial smooth solution $\Theta$
to the inviscid quasi-geostrofic equations~\eqref{eq:SQG transport}\eqref{eq:SQG transport} which has the form
\begin{equation}\label{eq:thrm}
\Theta (r, z, t) = \Theta (r, z - ct) = f\big(
\Psi (r, z - ct) - cr - k\big),
\end{equation}
for all $(r, z) \in \HH$ and that satisfies the symmetries
$\Theta (r, z) = -\Theta (-r, z) = \Theta (r, -z),$
for all $(r, z) \in \RR^2$. Moreover, The restriction of $\Theta$ to $\HH$ is non-negative, compactly supported and
non-increasing relatively to the variable $|z|$.
\end{theorem}

\section{Strategy of proof and main lemmas}
We regroup in this section the main Lemmas involved in the proof of Theorem~\ref{thrm:} and how they
follow one another. The detailed proof of these different lemmas are provided in Section~\ref{sec:proofs}.

\subsection{Properties of the Nehari Manifold and minimizing sequences}\label{sec:2_1}
We are interested in the minimization problem
\begin{equation}
\alpha := \inf\,\{E (\Psi) : \Psi \in \cN\} . \label{eq:minimization problem}
\end{equation}
Since the function $f$ worth $0$ on $\RR_-$ then a given function $\Psi$ cannot belong to $\cN$ if $\Psi\leq 0$ on $\HH$
Indeed, this would imply that
\begin{equation}
\int_{\HH}f (\Psi - cr - k) \Psi = 0
\end{equation}
and then $\|\Psi\|_{X^s}$.
The only function on $X^s_{sym}$ such that this quantity worth 0 is the null function which has been
excluded from the definition of the Nehari manifold). We have the following description of the
Nehari manifold. 
\begin{lemma}\label{lem:Nehari properties}
The set $\cN$ is a $\cC^1$ non-empty sub-manifold of $X^s_{sym}.$ For every $\Psi \in X^s_{sym}$ such
that $\cL^2( \supp(\Psi_+) \cap \HH)$ is non zero\footnote{The notation $\cL^d$ refers to the $d$-dimensional Lebesgue measure, ie. the lebesgue measure on $\RR^d$. The function $\Psi_+:=\max\{\Psi,0\}$ is the positive part of $\Psi$.}, there exist a unique $t_\Psi > 0$ such that $t_\Psi\Psi \in\cN$. The value of this $t_\Psi$ is characterized by
\begin{equation}\label{eq:characterization of t_psi} E (t_\Psi\Psi) = \max\{E (t\Psi) : t > 0\} .
\end{equation}
Moreover, any local minimizer of $E$ on $N$ is a smooth non-trivial solution of~\eqref{eq:frac lap semi-lin}. We also have that
\begin{equation}
\beta := \inf \big\{\|\Psi\|^2_{X^s}\;:\;\Psi\in\cN\big\}> 0. \label{eq:beta positive}
\end{equation}
and for every $\Psi \in\cN$ ,
\begin{equation}
\|\Psi\|^2_{X^s} \leq\bigg(1+\frac{1}{\mu}\bigg)E (\Psi). \label{eq:alpha positive} 
\end{equation}
\end{lemma}
Remark that this last assertion implies that $\alpha$ is positive. This proposition also implies that
any minimizing sequence of $E$ on $\cN$ is a bounded sequence.
\begin{definition}[Polarization]
We now define the \emph{polarization} of a function $\Psi \in X^s$ by
\begin{equation}
\forall \;x=(r,z)\in\RR^2,\qquad\Psi^\dag(X) :=\left\{\begin{array}{ll}
\max\big\{\Psi(x), \Psi\big(\sigma (x)\big)\big\}&\quad\text{if }r > 0,\\
\min\big\{\Psi(x), \Psi(\sigma (x)\big)\big\}&\quad\text{if }r < 0,
\end{array}\right.
\end{equation}
\end{definition}
where $\sigma$ denotes the linear map $x=(r, z) \in\RR^2\mapsto(-r,z)$.
In the particular case $\Psi\in X^s_{sym}$, we obtain $\Psi^\dag_{|\HH} \geq 0$ and $\Psi^\dag_{|\HH^c} \leq 0$. For more details about polarization, see for instance~\cite{Smets_Willem_2003}. 

\begin{lemma}[Polarization inequality]\label{lem:polarization inequality} 
For all $\Psi \in \cN$ , 
\begin{equation}
E(t_{\Psi^\dag}\Psi^\dag)\leq E (\Psi)
\end{equation}
and this inequality is strict when $\Psi \neq \Psi^\dag$.
\end{lemma}
Denote with a $\dag$ the image of a given set by the polarization. This lemma tells that if $(\Psi_n)$ is a minimizing sequence for $E$ on $\cN$ then so is $t_{\Psi^\dag}\Psi^\dag$ because by definition of $\Psi\mapsto t_\Psi$ the function $t_{\Psi^\dag}\Psi^\dag$ belongs to $\cN$.
Thus, the minimizer, if it exists, belongs to $\cN^\dag$.
It is then possible to restrict the investigations to $X^{s,\dag}_{sym}$.

\begin{definition}[Steiner rearrangement] We define the \emph{Steiner rearrangement} of $\Psi\in X^{s,\dag}_{sym}$, noted $\Psi^\sharp$, as being the function
of $X^{s,\dag}_{sym}$ which super-level sets on $\HH$ are given for all $\nu > 0$ by
\begin{equation}
\{\Psi^\sharp\geq \nu\} :=\bigcup_{r\in\RR+}\{r\}\times\Big[-\frac{\zeta_\Psi(r)}{2},+\frac{\zeta_\Psi(r)}{2}\Big]
\end{equation}
with
\begin{equation}
\zeta_\Psi(r) := L^1\{z \in\RR\,:\,\Psi(r, z) \geq \nu\}.
\end{equation}
We extend this definition on $\HH^c$ by symmetry to ensure that $\Psi^\sharp\in X^{s,\dag}_{sym}$.
\end{definition}

\begin{lemma}[Steiner inequality] \label{lem:Steiner inequality}
For all $\Psi \in\cN^\dag$, 
\begin{equation}
E(t_{\Psi^\sharp}\Psi^{\sharp})\leq E(\Psi)
\end{equation}
and the equality holds if and only if $\Psi=\Psi^{\sharp}$ up to a translation on the $z$ axis.
\end{lemma}
Then, if $(\Psi_n)$ is a minimizing sequence for $E$ on $\cN^\dag$ then so is $t_{\Psi^\sharp}\Psi^\sharp$. Thus, similarly as before it is possible to restrict the investigations to $X^{s,\sharp}_{sym}$.

\subsection{Existence of the solution of the minimizaing problem}\label{sec:2_2}
Let $(\Psi_n)\in\cN^\sharp$ a minimizing sequence. We already know that such a sequence in bounded as a consequence of Lemma~\eqref{lem:Nehari properties}.
To start with, we establish the following compactness result.

\begin{lemma}[compactness]\label{lem:compactness}
Let $c$ and $k$ be positive. Define the map
\begin{equation}
T : \Psi \in X^s
\longmapsto\left\{\begin{array}{ll}
(\Psi - cr - k)_+&\quad \text{on } \HH,\\
-(\Psi - cr + k)_-&\quad \text{on } \HH^c.
\end{array}\right.
\end{equation}
Then $T$ maps $X^s_{sym}$ into himself and maps bounded sets into bounded sets. Moreover, the map $T \circ \sharp \circ \dag$ is a compact map from $X^s_{sym}$ into $L^p_{sym}(\RR^2)$, with $1 \leq p <\frac{2}{1-s}$.
\end{lemma}
Up to an extraction we can suppose that the minimizing sequence $\Psi_n\to\Psi^\star$ weakly in $X^{s,\sharp}_{sym}$ and that $(\Psi_n -cr -k)_+ \to(\Psi^\sharp - cr - k)_+$ strongly in $L^p(\HH)$ for all $p <\frac{2}{1-s}$.

\begin{lemma}[convergence]\label{lem:convergence}
The convergence of $\Psi_n$ towards $\Psi^\star$ in $X^{s,\sharp}$ is a strong convergence.
\end{lemma}

This implies that $\Psi^\star$ is solution to the studied minimization problem.

\subsection{Properties of the solution}\label{sec:2_3}
We finally define $\Theta^\star$ from $\Psi^\star$ according to formula~\eqref{eq:ansatz}. Since $T(\Psi^\star) \in L^p(\RR^2)$ for all $p\in [1,\frac{2}{1-s}]$
then $\Theta^\star \in L^q$ for all $q \in [1,\frac{2}{\nu(1-s)}]$ as a consequence of~\eqref{eq:hypothesis 2}. We have the following regularity result

\begin{lemma}[regularity]\label{lem:regularity} The functions $\Psi^\star$ and $\Theta^\star$ are $\cC^\infty$. \end{lemma}
We can also establish a result on the
decay of $\Psi^\star$ at infinity.

\begin{lemma}[decay estimate]\label{lem:decay}
There exists a constant $C>0$ such that for all $r\in\RR^2$,
\begin{equation}
|\Psi^\star(x)| \leq \frac{C}{1 + |x|^{2(1-s)}},
\end{equation}
\end{lemma}
With the positive cut-off level$ k > 0$ appearing in the definition of $T$, this proposition implies in
particular that $\Theta^\star$ is compactly supported.

\section{Proofs of the lemmas}\label{sec:proofs}
\subsection{Proofs of the lemmas of section~\ref{sec:2_1}}
\subsubsection{Proof of Lemma~\ref{lem:Nehari properties}}
Let $\Psi \in X^s_{sym}$ with $\cL^2\big(\supp(\Psi_+) \cap \HH\big)\neq 0$. For any $t > 0$, we define
\begin{equation}
g (t) := \frac{E'(t\Psi) (t\Psi)}{t^2}=\frac{1}{2}\|\Psi\|^2_{X_s}-2t\int_{\HH}f(t\Psi_+ - cr - k) \Psi_+.
\end{equation}
We observe that the integral above can we rewritten
\begin{equation}
g(t) = \frac{1}{2}\|\Psi\|^2_{X^s} - 2\int_{\HH}\frac{f(t\Psi_+(r, z) - cr - k)}{t\Psi_+(r, z) - cr - k + (cr + k)}\;(\Psi+)^2
 (r, z)\, dr\, dz. \label{eq:rewrite the studied integral}
\end{equation}
Since we have $\cL^2\big(\supp(\Psi_+) \cap \HH\big)\neq 0$, our remark on the variations $\xi \mapsto f(\xi)/(\xi + \xi_0),$ consequence of \eqref{eq:hypothesis 3}, indicates that $t\mapsto g(t)$ is decreasing and $g(t)\to-\infty$ as $t\to+\infty$.
Indeed, one have to apply this property of $f$ to \eqref{eq:rewrite the studied integral} with $\xi = t\Psi+(r, z) - cr - k$ and $\xi_0 = cr + k$
and then integrate on $\HH$ against the non-negative weight $(\Psi_+)^2$.
We use Hypothesis \eqref{eq:hypothesis 2} to write on $\HH$
\begin{equation}
0 \leq\frac{1}{t}f (t\Psi_+ - cr - k) \Psi_+ \leq\frac{1}{t}f (t\Psi_+) \Psi_+ \leq Ct^{\nu-1}\big(\Psi_+\big)^{\nu+1}.
\end{equation}
Since $\nu>1$, 
\begin{equation}
g(t)\to\frac{1}{2}\|\Psi\|^2_{X^2} > 0.
\end{equation}
Since $f$ is smooth then $g$ is continuous, and then the function $g$ admits a unique root on $\RR_+^\ast$. 
The characterization \eqref{eq:characterization of t_psi} comes from the fact that
\begin{equation}
tg(t) = \frac{d}{dt}E (t\Psi).
\end{equation}
The estimate \eqref{eq:alpha positive} is obtained, for $\Psi \in \cN$, as follows
\begin{equation}\begin{split}
&E (\Psi) = E (\Psi) -\frac{1}{\mu + 1}E'(\Psi)(\Psi)\\
&=\frac{\mu}{2(\mu + 1)} \|\Psi\|^2_{X^s} +\frac{2}{\mu + 1} \int_{\HH}\Big[f (\Psi - cr - k) \Psi(r,z) - (\mu + 1)F (\Psi - cr - k)\Big]\,dr\,dz\\&\geq\frac{\mu}{2(\mu + 1)} \|\Psi\|^2_{X^s},
\end{split}
\end{equation}
The last inequality comes from the integration on $[0, r]$ of hypothesis \eqref{eq:hypothesis 3} that gives
$\mu F(t) \leq t f(t) - F(t).$
The fact that $\beta$ is not zero is obtained using~\eqref{eq:hypothesis 2} and the Sobolev embedding
\begin{equation}
\|\Psi\|^2_{X^s}=\int_{\HH}f (\Psi_+ - cr - k) \Psi_+\leq 4K\int_{\HH}(\Psi+)^\frac{2}{1-s}=4K \|\Psi_+\|^\frac{2}{1-s}_{L^\frac{2}{1-s}}\leq C \|\Psi\|^\frac{2}{1-s}_{X^s}.
\end{equation}
Concerning the regularity of $\cN$ , it is a consequence of the implicit functions theorem applied to
$\Xi : (s, \Psi) \mapsto E'(s\Psi) (\Psi)$
defined on the open set $\RR_+^\ast\times X^s_{sym}\setminus\{0\}.$
 The hypothesis of the theorem are verified because
for $\Psi \in \cN$ we have:
\begin{equation}
\partial_1\Xi (1, \Psi) = t^2_{\Psi} g'(t\Psi) < 0.
\end{equation}
It remains to prove that any minimizer of $E$ on $\cN$ is a critical point for $E$ defined on the whole
space. We first remark that a minimizer of $E$ on $N$ is a minimizer of $\Psi \mapsto E (t\Psi\Psi)$ on $X^s_{sym}$.
Then, using the definition of the Nehari manifold and the fact that we have $\Psi \in \cN$ implies $t_\Psi = 1$,
we conclude
\begin{equation}
\forall  h \in X^s_{sym}, E'(\Psi) (h) = E'(t\Psi)(\Psi) [t'_\Psi (h) \Psi + t_\Psi h] = 0.
\end{equation}\qed

\subsection{Proof of Lemma~\ref{lem:polarization inequality}}
We first recall that $\Psi \in \cN$ implies that
$\cL^2\big(\supp(\Psi_+) \cap \HH\big)\neq 0.$ Using the characterization \eqref{eq:characterization of t_psi} we get $E (\Psi) \geq E (t_{\Psi^\dag}\Psi).$
Using the fact that $\Psi(r, z) = -\Psi(-r, z),$ we conclude that here the polarization consists in switching the two values of $\Psi(r, z)$ and $\Psi(-r, z)$ if and only if we have $\Psi(r, z) \leq 0 \leq \Psi(-r, z).$ Therefore
since $F$ worth $0$ on $\RR_-$ and is positive on $\RR_+$, we obtain
\begin{equation}
V (t_{\Psi^\dag}\Psi) \leq V (t_{\Psi^\dag}\Psi^\dag
).
\end{equation}
To finish the proof of this lemma, we have to establish 
\begin{equation}
\|\Psi^\dag\|_{X^s} \leq \|\Psi\|_{X^s}
\end{equation}
and that this inequality is strict if and only if $\Psi^\dag \neq \Psi.$ Actually the fact that the polarization
decreases the $\dot{W}^{s,p}(\RR^d)$ half-norms \eqref{eq:gagliardo half-norms} is a general result so that we can establish it in the general
case. By definition of the principal values, we have
\begin{equation}
\iint_{|x-y|\geq\varepsilon}\frac{|\Psi(x) - \Psi(y)|^p}{|x - y|^{d+sp}}\,dx\,dy \longrightarrow |u|^p_{W^{s,p}} qquad\text{as }\varepsilon\to 0. \label{eq:principal values convergence} 
\end{equation}
We then establish the inequality for any fixed $\varepsilon > 0.$ First, the integral is split as follows,
\begin{equation}\begin{split}
&\iint_{|x-y|\geq\varepsilon}\frac{|\Psi(x) - \Psi(y)|^p}{|x - y|^{d+sp}}\,dx\, dy\\
&=\iint_{\HH^2\setminus\{|x-y|<\varepsilon\}}
\bigg(\frac{1}{|x-y|^{d+sp}}\Big(|\Psi(x)-\Psi(y)|^p+|\Psi\circ\sigma(x)-\Psi\circ\sigma(y)|^p\Big)\\
&\qquad+\frac{1}{|x-\sigma(y)|^{d+sp}}\Big(|\Psi(x)-\Psi\circ\sigma(y)|^p+|\Psi\circ\sigma(x)-\Psi(y)|^p\Big)\bigg)\,dx\,dy
\end{split}
\label{eq:split integral into 4} 
\end{equation}
Let $x, y \in \HH$. Observe that
\begin{equation}
|x - y|^{d+sp} < |x - \sigma(y)|^{d+sp}. \label{eq:estimate the distances} 
\end{equation}
\textbf{Case 1:} $\Psi(x) \geq \Psi \circ \sigma(x)$ and $\Psi(y) \geq \Psi \circ \sigma(y).$
In this case, with the definition of the polarization, $\Psi(x) = \Psi^\dag(x)$ and $\Psi(y) = \Psi\dag
(y)$.
Then when we integrate on the couples $(x, y)$ that belongs to Case 1, the associated term in the
integral \eqref{eq:split integral into 4} is not modified by the polarization.

\textbf{Case 2:} $\Psi(x) \geq \Psi \circ \sigma(x)$ and $\Psi(y) < \Psi \circ \sigma(y).$
By computing its derivative, we obtain that the function
\begin{equation}
u_{\beta,s } : \alpha \in \RR \mapsto |\alpha + \beta|^p - |\alpha + s |^p
\end{equation}
is non-decreasing when $\beta > s $. Indeed we have (with $p \geq 1$)
\begin{equation}
u'_{\beta,s }(\alpha) = p(\alpha + \beta)|\alpha + \beta|^{p-2} - p(\alpha + s )|\alpha + s |^{p-2}
(31)
\end{equation}
which is non-negative because $y \mapsto y|y|^{p-2}$
is an non-decreasing function. We now use this property
of $u_{\beta,s }$ with $\alpha_1 := \Psi(x) \geq \Psi \circ \sigma(x) =: \alpha_2$ and with $\beta := -u(y) > s  := -u \circ \sigma(y)$. We obtain
\begin{equation}
|\Psi(x)-\Psi(y)|^p+|\Psi\circ\sigma(x)-\Psi\circ\sigma(y)|^p>|\Psi\circ\sigma(x) - \Psi(y)|^p+|\Psi(x) - \Psi \circ \sigma(y)|^p.
\end{equation}
If we combine this with \eqref{eq:estimate the distances} we get
\begin{equation}\begin{split}
&\frac{1}{|x - y|^{d+sp}}\Big(|\Psi(x) - \Psi(y)|^p+|\Psi \circ \sigma(x) - \Psi \circ \sigma(y)|^p\Big)\\
&\qquad+\frac{1}{|x - \sigma(y)|^{d+sp}}\Big(|\Psi \circ \sigma(x) - \Psi(y)|^p+|\Psi(x) - \Psi \circ \sigma(y)|^p\Big)\\
&>\frac{1}{|x - y|^{d+sp}}\Big(|\Psi \circ \sigma(x) - \Psi(y)|^p+|\Psi(x) - \Psi \circ \sigma(y)|^p\Big)\\
&\qquad+\frac{1}{|x - \sigma(y)|^{d+sp}}\Big(|\Psi(x) - \Psi(y)|^p+\Psi \circ \sigma(x) - \Psi \circ \sigma(y)|^p\Big)\\
&=\frac{1}{|x - y|^{d+sp}}\Big(|\Psi^\dag(x) - \Psi^\dag
(y)|^p+|\Psi^\dag\circ \sigma(x) - \Psi^\dag\circ \sigma(y)|^p\Big)\\
&\qquad+\frac{1}{|x - \sigma(y)|^{d+sp}}\Big(|\Psi^\dag\circ \sigma(x) - \Psi^\dag(y)|^p+\Psi^\dag(x) - \Psi^\dag
\circ \sigma(y)|^p\Big).\end{split}
\end{equation}

\textbf{Case 3:} $\Psi(x) < \Psi \circ \sigma(x)$ and $\Psi(y) < \Psi \circ \sigma(y).$
In this case we have both $\Psi(x)$ and $\Psi(y)$ that are swapped with respectively $\Psi\circ\sigma(x)$ and $\Psi\circ\sigma(y).$
Then this case is the same as case 1 and the term associated to Case 3 in the integral \eqref{eq:split integral into 4} is not
modified by the polarization.

\textbf{Case 4:} $\Psi(x) < \Psi \circ \sigma(x)$ and $\Psi(y) \geq \Psi \circ \sigma(y).$
This case is the same as Case 2.

Gathering these four cases we obtain that for any $\varepsilon > 0,$
\begin{equation}
\iint_{|x-y|\geq\varepsilon}\frac{|\Psi(x) - \Psi(y)|^p}{|x - y|^{d+sp}}\,dx\,dy \geq\iint_{|x-y|\geq\varepsilon}\frac{|\Psi^\dag(x) - \Psi^\dag(y)|^p}{|x - y|^{d+sp}}\,dx\,dy. \label{eq:inequality with epsilon} 
\end{equation}
Concerning the cases of equality, we obtained from Cases 2 and 4 that if
\begin{equation}
\cL^2\bigg( \Big\{(x, y) \in \HH^2:\Psi(x) = \Psi^\dag(x)\quad\mathrm{and}\quad \Psi(y) \neq \Psi^\dag(y)\Big\}\cap {|x - y| \geq \varepsilon}\bigg)> 0
\end{equation}
then the inequality\eqref{eq:inequality with epsilon}  is actually strict. We now observe that the above set is of measure zero for
every $\varepsilon > 0$ if and only if we have either $\Psi = \Psi^\dag$ or $\Psi = \Psi^\dag \circ \sigma$. But this last case is not possible
when $\Psi \in \cN$ and then the only case of equality in our case is $\Psi = \Psi^\dag$.\qed

\subsection{Proof of Lemma~\ref{lem:Steiner inequality}}
Arguing similarly as the previous proof, we only have to prove that
\begin{equation}
\forall \; \Psi \in X^{s,\dag}_{sym},\quad E(\Psi^\sharp) \leq E(\Psi).
\end{equation}
Since the Steiner rearrangement only involves rearrangements of the super-level sets perpendicularly to the $r$-axis, we get
\begin{equation}
V(\Psi^\sharp) = V(\Psi).
\end{equation}
To conclude we have to establish that.
\begin{equation}
\|\Psi^\sharp\|_{X^s} \leq \|\Psi\|_{X^s}.
\end{equation}
To start with, we suppose that $\Psi$ is smooth and compactly supported. In this case
\begin{equation}
\int_{\RR^2}\int_{\RR^2}\frac{|\Psi(x)-\Psi(y)|^2}{(|x - y|^2 + \varepsilon^2)^{1+s}}\,dx\,dy \;\longrightarrow\; \|\Psi\|_{X^s},
\end{equation}
as $\varepsilon\to0^+$. Since the considered functions are $\cC^\infty$, is is possible to develop the square above and write
\begin{equation}\begin{split}
&\qquad\qquad\qquad\int_{\RR^2}\int_{\RR^2}\frac{|\Psi (x) - \Psi (y)|^2}{(|x - y|^2 + \varepsilon^2)^{1+s}}\,dx\, dy\\
&= 2 \int_{\RR^2}\int_{\RR^2}\frac{\Psi (x)^2}{(|x - y|^2 + \varepsilon^2)^{1+s}}\,dx\,dy - 2\int_{\RR^2}\int_{\RR^2}\frac{\Psi (x) \Psi (y)}{(|x - y|^2 + \varepsilon^2)^{1+s}}\,dx\,dy. \end{split}\label{estim_hs} 
\end{equation}
The first integral in the right-hand side of the above inequality is not modified by rearrangement
of the super-level sets of the function $\Psi$. Concerning the second integral, using the fact that $\Psi(r, z) = -\Psi(-r, z)$ we get
\begin{equation}\begin{split}
&\int^{+\infty}_{-\infty}\int^{+\infty}_{-\infty}\int^{+\infty}_{-\infty}\int^{+\infty}_{-\infty}\frac{\Psi(r_x, z_x)\Psi(r_y, z_y)}{((r_x - r_y)^2 + (z_x - z_y)^2 + \varepsilon^2)^{1+s}}\,dz_x\,dz_y\,dr_x\,dr_y\\
&= 2 \int^{+\infty}_{0}\int^{+\infty}_{0}\int^{+\infty}_{-\infty}\int^{+\infty}_{-\infty}
\frac{\Psi(r_x, z_x)\Psi(r_y, z_y)}{((r_x - r_y)^2 + (z_x - z_y)^2 + \varepsilon^2)^{1+s}}\,dz_x\,dz_y\,dr_x\,dr_y \\
&- 2\int^{+\infty}_{0}\int^{+\infty}_{0}\int^{+\infty}_{-\infty}\int^{+\infty}_{-\infty}\frac{\Psi(r_x, z_x)\Psi(r_y, z_y)}{((r_x + r_y)^2 + (z_x - z_y)^2 + \varepsilon^2)^{1+s}}\,dz_x\,dz_y\,dr_x\,dr_y 
\end{split}\label{il_faut_estimer}\end{equation}
But the function
\begin{equation}
\Upsilon_{r_x,r_y}: u\;\longmapsto\;\frac{1}{\big((r_x - r_y)^2 + u^2 + \varepsilon^2\big)^{1+s}}-\frac{1}{\big((r_x + r_y)^2 + u^2 + \varepsilon^2\big)^{1+s}}
\end{equation}
is non-negative and radially decreasing on $\RR$. Moreover, for $r_x, r_y \geq 0$ the functions $z_x \mapsto \Psi(r_x, z_x)$
and $z_y \mapsto \Psi(r_y, z_y)$ are both non-negative $\RR$. Thus, using the Riesz rearrangement inequality, we
obtain
\begin{equation}\begin{split}
&\int^{+\infty}_{-\infty}\int^{+\infty}_{-\infty}\Psi(r_x, z_x)\Psi(r_y, z_y)\Upsilon_{r_x,r_y}(z_x - z_y) \,dz_x\, dz_y \\
&\qquad\leq\int^{+\infty}_{-\infty}\int^{+\infty}_{-\infty}\Psi^\sharp(r_x, z_x)\Psi^\sharp(r_y, z_y)\Upsilon_{r_x,r_y}(z_x - z_y) \,dz_x\, dz_y.
\end{split}
\end{equation}
We now inject this inequality back into \eqref{il_faut_estimer} and get
\begin{equation}
\int_{\RR^2}\int_{\RR^2}\frac{\Psi (x) \Psi (y)}{(|x - y|^2 + \varepsilon^2)^{1+s}}\,dx\,dy \leq\int_{\RR^2}\int_{\RR^2}\frac{\Psi^\natural (x) \Psi^\natural (y)}{(|x - y|^2 + \varepsilon^2)^{1+s}}\,dx\,dy.
\end{equation}
We use this estimate in \eqref{estim_hs}, we take the limit and we conclude by density of the smooth compactly supported functions. \qed

Remark: It was not possible to use directly the Riesz rearrangement inequality to the second
integral appearing in \eqref{estim_hs} because this inequality in only true for non-negative functions.\qed

\subsection{Proofs of the lemmas of section~\ref{sec:2_2}}
\subsubsection{Proof of Lemma~\ref{lem:compactness}}
~\newline
\indent\textbf{\mathversion{bold}$\bullet\quad$Step 1 : $T$ maps $X^{s,\dag}_{sym}$ into itself and maps bounded sets into bounded sets.}\newline
First, if $\Psi$ satisfies the symmetry property then so does $T(\Psi)$. 
Define the set\footnote{The adherence of this set is the support of the function $\Theta$. This corresponds physically speaking to the vorticity zone.}
\begin{equation}
\Omega(\Psi) := \{(r, z) \in \HH: T(\Psi)(r, z) > 0\},
\end{equation}
where $T(\Psi)(r, z):=(\Psi(r,z)-cr-k)_+$. By definition of $T$ and of $\Omega$
\begin{equation}
\cL^2(\Omega) = \int_\Omega1 \leq\int_\Omega\bigg(\frac{\Psi(r,z)}{cr + k}\bigg)^\frac{2}{1-s}dr\,dz\leq\int_\Omega\bigg(\frac{\Psi(r,z)}{k}\bigg)^\frac{2}{1-s}dr\,dz\leq
\frac{1}{k^\frac{2}{1-s}}\int_\HH\Psi^\frac{2}{1-s}. \label{estimOmega}
\end{equation}
Using a Sobolev inequality above leads to
\begin{equation}
\cL^2(\Omega) \leq \frac{C_s}{k^\frac{2}{1-s}}\|\Psi\|_{X^s}^\frac{2}{1-s}.
\label{psigam}
\end{equation}
The computation of the double integral defining the $\dot{H}^s$ half-norm~\ref{eq:gagliardo half-norms} is done separating the integrals on $\RR^2$ on two between $\Omega$ and $\Omega^c$. On $\Omega^c$ the quantity $\Psi(r,z)-cr-k$ is non-positive and then $T(\Psi)(r,z)=0$. Therefore,
\begin{equation}\label{eq:estim}
\int_{\Omega^c}\int_{\Omega^c}\frac{|T(\Psi)(x) - T(\Psi)(y)|^2}{|x - y|^2(1+s)}\,dx\,dy = 0.
\end{equation}
Concerning the integral on $\Omega\times\Omega$, using the notation $x=(r_x,z_x)$ and $y=(r_y,z_y)$,
\begin{equation}
\begin{split}
\int_\Omega\int_\Omega\frac{|T(\Psi)(x) - T(\Psi)(y)|^2}{|x - y|^{2(1+s)}}\,dx\,dy&=\int_\Omega\int_\Omega\frac{|(\Psi(x) - cr_x - k) - (\Psi(y) - cr_y - k)|^2}{|x - y|^{2(1+s)}}\,dr\,dz\\&\leq\int_{\Omega}\int_\Omega\frac{|\Psi(x) - \Psi(y)|^2 + c^2|r_x - r_y|^2}{|x - y|^{2(1+s)}}
dr dz.
\end{split}
\end{equation}
Denote with an $\ast$ the radially decreasing rearrangement and $R_\Omega > 0$ the radius such that
\begin{equation}\cL^2(\Omega) = \cL^2 ( \cB (0, R_\Omega) ).
\end{equation}
To simplify the notations we simply note this ball $B(\Omega)$. By the Riesz rearrangement inequality
we have
\begin{equation}\begin{split}
\int_\Omega\int_\Omega\frac{|r_x - r_y|^2}{|x - y|^{2(1+s)}}\,dx\,dy&\leq\int_\Omega\int_\Omega\frac{dx\, dy}{|x - y|^{2s}}\leq\int_{\cB(\Omega)}\int_{\cB(\Omega)}\frac{dx\, dy}{|x - y|^{2s}}\\&\leq
\int_{\cB(\Omega)}\int_{\cB(\Omega)}\frac{dx}{|x|^{2s}}dy =\frac{\pi^s}{1-s}\cL^2(\Omega)^{2-s}.
\end{split}
\end{equation}
Using now \eqref{psigam} we get
\begin{equation}
\int_\Omega\int_\Omega\frac{|r_x - r_y|^2}{|r - z|^{2(1+s)}}\,dr\,dz \leq C_s\,\|\Psi\|_{X^s}^{2\frac{2-s}{1-s}}. 
\end{equation}
Concerning the last term,
\begin{equation}
\int_\Omega\int_{\Omega^c}\frac{|T(\Psi)(x) - T(\Psi)(y)|^2}{|x - y|^{2(1+s)}}\,dx\,dy =\int_\Omega|T(\Psi)(x)|^2\int_{\Omega^c}\frac{dy}{|x - y|^{2(1+s)}} dx. \label{dernier}
\end{equation}
For all $x \in \Omega$ we define $\displaystyle\Lambda(x) := \frac{\Psi(x) - cr_x - k}{2c}$ and $\cO_x := \{y \in \Omega^c : r_y \geq r_x + \Lambda(x)\}$. Then, 
\begin{equation}
\int_{\cO_x}\frac{dy}{|x - y|^{2(1+s)}}\leq
\int_{\cO_x}\frac{dy}{|r_x - r_y|^{2(1+s)}}\leq\int_{\cB(x,\Lambda(x))^c}\frac{dy}{|x - y|^{2(1+s)}}=\frac{\pi^s}{\Lambda(x)^{2s}}. \label{unsurlambda}
\end{equation}
Using the factq that $x \in \Omega$, that $y \in \Omega^c$ and that $y\in \cO_x$ in this order, we obtain
\begin{equation}\begin{split}
&\quad0 \leq T(\Psi)(x) = \Psi(x) - cr_x - k\leq \Psi(x) - \Psi(y) + c(r_x - r_y)\\&\leq \Psi(x) - \Psi(y) + \Psi(x) - \frac{cr_x - k}{2}= \Psi(x) - \Psi(y) + \frac{1}{2}T(\Psi)(x).\end{split}
\end{equation}
Therefore
\begin{equation}
|T(\Psi)(x)| \leq 2|\Psi(x) - \Psi(y)| \label{majorT} 
\end{equation}
Combing \eqref{dernier}, \eqref{unsurlambda} and \eqref{majorT} leads to
\begin{equation}\begin{split}
&\qquad\qquad\int_\Omega\int_{\Omega^c}\frac{|T(\Psi)(x) - T(\Psi)(y)|^2}{|x - y|^{2(1+s)}}dx dy\\&\leq 4\int_\Omega\int_{\Omega^c\setminus\cO_x}\frac{|\Psi(x) - \Psi(y)|^2}{|x - y|^{2(1+s)}}\,dx\, dy +\int_\Omega|T(\Psi)(x)|^2\frac{\pi}{s\Lambda(x)^{2s}}dx\\
&\leq 4\int_\Omega\int_{\Omega^c}\frac{|\Psi(x) - \Psi(y)|^2}{|x - y|
^{2(1+s)}}\,dx\, dy +\frac{\pi}{s}\int_\Omega|T(\Psi)(x)|^{2(1-s)}dx.
\end{split}
\end{equation}
Now, to estimate the last term of the above inequality, we use the fact that $T(\Psi) \leq \Psi$ and then the Hölder inequality gives
\begin{equation}
\int_\Omega|T(\Psi)(x)|^{2(1-s)}dx \leq\int_\Omega|\Psi(x)|^{2(1-s)}dx \leq \cL^2(\Omega)^{s(2-s)}\,\|\Psi\|^{2(1-s)}_{L^\frac{2}{1-s}}.
\end{equation}
We continue this estimate using \eqref{psigam} and a Sobolev embedding,
\begin{equation}
\leq C\|\Psi\|^{2s(1+\frac{1}{1-s} )}_{X^s} \|\Psi\|^{2(1-s)}_{L^\frac{2}{1-s}}\leq C\|\Psi\|^\frac{2}{1-s}_{X^s} . 
\end{equation}
Thus, gathering all these estimates we obtain.
\begin{equation}
\|T(\Psi)\|^2_{X^s}:=\int_\Omega\int_{\Omega^c}\frac{|T(\Psi)(x) - T(\Psi)(y)|^2}{|x - y|^{2(1+s)}}\,dx\, dy \leq C\|\Psi\|^2_{X^s}\Big(1 + \|\Psi\|^\frac{2s}{1-s}_{X^s}\Big).
\end{equation}
Therefore, $T$ does map $X^{s,\dag}_{sym}$ into itself and maps bounded subsets of $X^{s,\dag}_{sym}$ into bounded subsets.\newline

\textbf{\mathversion{bold}$\bullet\quad$Step 2 : $T\circ\sharp\circ\dag$ defined on $X^s_{sym}$ is a compact operator for the $L^p$ topology.}
Set the convention that $\{|z| \geq R\}$ designates the set $\{(r, z) \in \HH : |z| \geq R\}$. Let $\kappa > 0$ and $R \geq 0$. For all $r \in \Omega^\sharp$ we define
\begin{equation}
\cU_x := \cB (x, \kappa) \cap \big(\Omega^\sharp\big)^c.
\end{equation}
Then,
\begin{equation}\begin{split}
\int_{\{|z|\geq R\}}&|T(\Psi^\sharp)(x)|^2dx =\int_{\{|z|\geq R\}}\frac{1}{\cL^2(\cU_x)}\int_{\cU_x}|T(\Psi^\sharp)(x)|^2\,dy\,dx\\
&\leq\int_{\{|z|\geq R\}}\frac{1}{\cL^2(\cU_x)}\int_{\cU_x}|T(\Psi^\sharp)(x) - T(\Psi^\sharp)(y)|^2\, dy\,dx\\
&\leq\int_{\{|z|\geq R\}}\frac{\kappa^{2(1+s)}}{\cL^2(\cU_x)}\int_{\cU_x}\frac{|T(\Psi^\sharp)(x) - T(\Psi^\sharp)(y)|^2}{|x - y|^{2(1+s)}}\,dy\,dx. 
\end{split}\label{enz1} 
\end{equation}
Denote by $P_\RR$ the projection on $\RR \times {0}$ (that is identified to $\RR$). As a consequence of the Steiner symmetrization, with \eqref{estimOmega},
\begin{equation}
2R\,\cL^1\bigg(P_\RR\Big((\Omega^\sharp)^c\cap \{|z| \geq R\}
\Big)\bigg) \leq \cL^2\big(\Omega^\sharp\big)\leq\frac{1}{k^\frac{2}{1-s}}\|\Psi^\sharp\|^\frac{2}{1-s}_{L^\frac{2}{1-s}}. \label{estimBande} 
\end{equation}

Since $|z_x| \geq R-\kappa$ then using again the Steiner symmetry of $\Omega^\sharp$, gives that $\cU_x$ contains the ball $B(r, \kappa)$ minus the rectangle centered at $x$, of width $$\cL^1\Big(P_\RR\big((\Omega^\sharp)^c\cap \{|z| \geq R - \kappa\}\big)\Big)$$
height $2\kappa$. Then, with \eqref{estimBande},
\begin{equation}
\cL^2 (\cU_x) \geq \pi\kappa^2 -\frac{\kappa}{(R - \kappa) k^\frac{2}{1-s}}\|\Psi^\sharp\|^\frac{2}{1-s}_{L^\frac{2}{1-s}}. \label{Ux} 
\end{equation}
The choice of $\kappa$ is free and then we choose to fix it equal to $C/R$ with 
$$C :=\frac{4}{\pi\, k^\frac{2}{1-s}}\|\Psi^\sharp\|^\frac{2}{1-s}_{L^\frac{2}{1-s}}.$$ Since
Choose now $R$ such that $R \geq\sqrt{2C}$. Then in this case the inequality \eqref{Ux} becomes
\begin{equation}
\cL^2(\cU_x) \geq\frac{\pi C^2}{2R^2}. \end{equation}
Combining the estimate above with \eqref{enz1}, leads to the following estimate
\begin{equation}
\int_{\{|z|\geq R\}}|T(\Psi^\sharp)(x)|^2dx \leq\bigg(\frac{4}{\pi R}\bigg)^{2s}\Bigg(\frac{\|\Psi^\sharp\|_{L^\frac{2}{1-s}}}{k}\Bigg)^\frac{4s}{1-s}\|T(\Psi^\sharp)\|^2_{X^s}  . \label{enz}
\end{equation}
On the other hand, using the Hölder inequality,
\begin{equation}\label{pif}\begin{split}
\int_{\{r\geq R\}}|T(\Psi^\sharp)|^2=\int_{\{r\geq R\}}|T(\Psi^\sharp)|^2\mathbbm{1}_{(\Omega^\sharp)^c}&\leq\bigg(\int_{\{r\geq R\}}|T(\Psi^\sharp)|^\frac{2}{1-s}\bigg)^{1-s}\bigg(\int_{\{r\geq R\}}\mathbbm{1}_{(\Omega^\sharp)^c}\bigg)^s\\
&=cL^2\Big((\Omega^\sharp)^c\cap\{r\geq R\}\Big)^s\|\Psi^\sharp\|^2_{L^\frac{2}{1-s}},
\end{split}
\end{equation}
where by convention $\{r \geq R\}$ designates the set $\{(r, z) \in\HH : r \geq R\}$. Moreover
\begin{equation}\label{paf}\begin{split}
&\qquad\cL^2\big((\Omega^\sharp)^c\cap \{r \geq R\}\big) = \int_{(\Omega^\sharp)^c\cap \{r \geq R\}}1\\
&\leq\int_{(\Omega^\sharp)^c\cap \{r \geq R\}}\bigg(\frac{\Psi^\sharp}{cr + k}\bigg)^\frac{2}{1-s}\leq\bigg(\frac{1}{cR + k}\bigg)^\frac{2}{1-s}\|\Psi^\sharp\|^\frac{2}{1-s}_{L^\frac{2}{1-s}}.
\end{split}
\end{equation}
Combining~\eqref{pif} and~\eqref{paf} leads to
\begin{equation}
\int_{\{r\geq R\}}|T(\Psi^\sharp)|^2\leq\bigg(\frac{1}{cR + k}\bigg)^\frac{2s}{1-s}\|\Psi^\sharp\|^\frac{2}{1-s}_{L^\frac{2}{1-s}}\label{enr}
\end{equation}

The two decay estimates \eqref{enz} et \eqref{enr} and the Rellich-Kondrachov compactness theorem (applied
at the local level) give the result.\qed

\subsubsection{Proof of Lemma~\ref{lem:convergence}}
It follows from the definition of $\cN$ and of Lemma~\ref{lem:Nehari properties} that
\begin{equation}
\int_\HH f(\Psi_n - cr - k)\Psi_n =\frac{1}{2}\int_\RR\Psi_n(-\Delta)^s\Psi_n =\frac{1}{2}\|\Psi_n\|^2_{X^s} \geq\frac{\beta}{2}>0.
\end{equation}
By the previous lemma, up to a subsequence when $n\to+\infty$,
\begin{equation}
\int_\HH f(\Psi_n - cr - k)\Psi_n \quad\longrightarrow\quad \int_\HH f(\Psi^\star - cr - k)\Psi\star \;\geq\frac{\beta}{2}.
\end{equation}
In particular $(\Psi^\star -cr -k)+ \not\equiv 0$ on $\HH$. 
By proposition~\ref{lem:Nehari properties}, there exists $t^\star > 0$ such that $t^\star\Psi^\star \in \cN$ .
With the characterization of $t^\star$ and since $\Psi_n \in \cN$,
\begin{equation}
E(\Psi_n) = E(t_{\Psi_n}\Psi_n) \geq E(t^\star\Psi_n).
\end{equation}
Thus,
\begin{equation}
\alpha = \lim\limits_{n\to+\infty}E(\Psi_n) \geq \liminf\limits_{n\to+\infty}E(t^\star\Psi_n) \geq E(t^\star\Psi^\star) \geq \alpha.
\end{equation}
Therefore all these inequalities are equalities and $\|\Psi_n\|^2_{X^s}\to|\Psi^\star\|^2_{X^s}$ . Since the space $X^s$ is
strictly convex, this gives that $\Psi_n$ converges towards $\Psi^\star$ stronly in $X^s$.\qed

\subsection{Proofs of the lemmas of section~\ref{sec:2_3}}
\subsubsection{Proof of Proposition~\ref{lem:regularity}}
We already know that $T(\Psi^\star )\in L^\frac{2}{1-s }(\RR2)$. Since the support of $T(\Psi\star )$ has finite measure, $T(\Psi^\star )\in L^1(\RR^2)$. Define $\Theta^\star$ from $\Psi^\star$ using formula~\eqref{eq:thrm}. Hypothesis~\eqref{eq:hypothesis 2} implies
\begin{equation}
\forall q\in\Big[1,\frac{2}{\nu (1-s)}\Big],\quad\Theta \star \in L^q(\RR^2).
\end{equation}
Define now the function $\widetilde{\Psi}$ given by the following representation formula,
\begin{equation}\label{eq:integral}
\widetilde{\Psi}(x) =K_s  \int_{\RR^2}\frac{\Theta^\star (y)}{|x-y|^{2(1-s)}}\,dy.
\end{equation} 
where $K_s$ is some renormalization constant.
It follows from the weighted inequalities for singular integrals~\cite[§5]{Stein_1993} that $\widetilde{\Psi}\in\dot{W}^{2s,q}(\RR^2)$, for all $q\in[1,\frac{2}{\nu (1-s)}].$ Moreover, by the Hardy-Littlewood-Sobolev convolution inequality,
$\widetilde{\Psi}\in L^q(\RR2),\forall q\in[\frac{1}{1-s},\frac{2}{\nu -s  (2 +\nu )}]$.
By standard interpolation, $\widetilde{\Psi}\in X^s_{sym}$.
Now, let $\varphi \in X^s_{sym}$ be a test function. Using the spectal properties of the Sobolev spaces~\cite{DiNezza_Palatucci_Valcinoci_2012} gives (up to multiplicative renormalization constants),
\begin{equation}\begin{split}
\left<\widetilde{\Psi},\varphi\right>_{X^s}&=\int_{\RR^2}|\xi|^{2s}\,\cF[\widetilde{\Psi}](\xi)\,\cF[\varphi](\xi)\,d\xi\\
&=\int_{\RR^2}|\xi|^{2s}\cF\Big[\Theta^\star\ast\frac{1}{|.|^{2(1-s}}\Big](\xi)\,\cF[\varphi](\xi)\,d\xi\\
&=\int_{\RR^2}\cF[\Theta^\star](\xi)\,\cF[\varphi](\xi)\,d\xi=\left<\Theta^\star,\varphi\right>_{L^2},
\end{split}\end{equation}
where $\cF[.]$ designates the Fourier transform.
Moreover, since $\Psi^\star$ is a critical point of $E$, then $\left<\Psi^\star ,\varphi \right>_{X^s} =\left<\Theta^\star ,\varphi \right>_{L^2}$, which implies $\widetilde{\Psi} = \Psi^\star$ .
The regularity known for $\Theta^\star$ allows to conclude that $\Psi^\star$ is bounded and uniformly continuous. Seen the definitions, to conclude that $\Psi^\star$ is smooth by a bootstrap argument there remain to study possible discontinuities on $r=0$. Nevertheless, it follows from the symmetry property of $\Psi^\star$ and its unifom continuity that $T(\Psi\star )$ worth $0$ at a distance uniformly positive from $r= 0$, meaning on a strip $]-\delta ,\delta [\times \RR_+.$ Therefore so is the case for $\Theta^\star$ and then the smoothness of $\Psi^\star$ is proved.\qed

\subsubsection{Proof of Proposition~\ref{lem:decay}}
Let $r \in R^2$ such that $|r| \geq 1$. We separate the integral~\eqref{eq:integral} into two,
\begin{equation}\label{eq:two integrals}
\Psi^\star(x) = K_s \int_{|x-y|\leq \frac{|x|}{2}}\frac{\Theta^\star(y)}{|x - y|^{2(1-s)}}\,dy +K_s \int_{|x-y|>\frac{|x|}{2}}\frac{\Theta^\star(y)}{|x - y|^{2(1-s)}}\,dy 
\end{equation}
Concerning the first integral, we choose $\eta \in ]1-s ,\frac{1}{s + 1}[$. This interval is non-empty and included in $]0, 1[$. We use the Hölder inequality and Hypothesis \eqref{eq:hypothesis 2} and then we are led to\begin{equation}\label{eq:holder used}
\int_{|x-y|\leq \frac{|x|}{2}}\frac{\Theta^\star(y)}{|x - y|^{2(1-s)}}\,dy\leq
 \bigg(\int_{|\zeta|\leq \frac{|x|}{2}}\frac{d\zeta}{|\zeta|^\frac{2(1-s)}{\eta}}\bigg)^\eta\bigg(\int_{|x-y|\leq\frac{|x|}{2}}|T (\Psi^\star)|^\frac{\nu}{1-\eta}\bigg)^{1-\eta}.
\end{equation}
Using again the estimates~\eqref{enz} et~\eqref{enr},
\begin{equation}
\int_{|x-y|\leq\frac{|x|}{2}}|T(\Psi^\star)|^2(y)dy \leq\frac{C}{(1 + |x|)^{-2s}}
. 
\end{equation}
Knowing that $\frac{\nu}{1 - \eta}\geq 2$ the above estimate used in~\eqref{eq:holder used} leads to (the constant that depends only on $\eta$)
\begin{equation}
\int_{|x-y|\leq \frac{|x|}{2}}\frac{\Theta^\star(y)}{|x - y|^{2(1-s)}}\,dy \leq C(\eta)\big(|x|^2 + 1\big)^{\eta(s+1)-1.} \label{estim1} 
\end{equation}
The second integral in~\eqref{eq:two integrals} can be estimated using directly the hypothesis on the function $f$,
\begin{equation}
\int_{|x-y|>\frac{|x|}{2}}\frac{\Theta^\star(y)}{|x - y|^{2(1-s)}}\,dy\leq\bigg(\frac{2}{|x|}\bigg)^{2(1-s)}\int_{\RR^2}|T(\Psi^\star)|^\nu \leq\frac{C}{|x|^{2(1-s)}}.\label{estim2}
\end{equation}
By choosing $\eta\in ]1 - s ,\frac{1}{s + 1}[$ such that $\eta \geq\frac{s}{s + 1}$, the estimates~\eqref{estim1} et~\eqref{estim2} give
\begin{equation}
\Psi^\star(x) \leq\frac{C}{1 + |x|^{2(1-s)}}.
\end{equation}

\bibliographystyle{plain}
\bibliography{bibliography}
\vspace{0.5cm}
The author acknowledges grants from the Agence nationale de la Recherche,  for project ``Ondes Dispersives 
Aléatoires" (ANR-18-CE40-0020-01).

\end{document}